\documentclass[a4paper,11pt,twoside]{article}
\usepackage{pst-fill,pst-grad}
\usepackage{textcomp} 
\usepackage[english]{babel}  
\usepackage[utf8]{inputenc}  
\usepackage[titletoc]{appendix}
\usepackage{titlesec}
\usepackage{graphicx} 
\usepackage{amsmath}
\usepackage{float}  
\usepackage{fancyhdr} 
\usepackage[matrix,arrow,curve]{xy} 
\usepackage{pstricks} 
\usepackage{amsmath,amsfonts,verbatim,afterpage,theorem,euscript,mathrsfs,amssymb}
\usepackage{amsfonts}
\usepackage{amssymb}
\usepackage{array} 
\usepackage{dsfont}
\usepackage[colorlinks=true,linkcolor=blue,citecolor=red]{hyperref}
\usepackage{authblk}
\usepackage{color}
 
\newtheorem{Proposition}{Proposition}[section]

\newtheorem{TheoremeP}{Theorem}
\newtheorem{CorollaireP}{Corollary}

\def \U{\vec{U}}
\def \W{\vec{W}}
\def \V{\vec{V}}

\def \Rt{\mathbb{R}^3}

\def \finpv{\hfill $\blacksquare$  \\ \newline }
\def \pv{{\bf{Proof.}}~} 

\def \ds{\displaystyle}

\title{\bf A short note on the Liouville problem for the steady-state Navier-Stokes equations}

\author[1]{ Oscar Jarr\'in\footnote{corresponding author: oscar.jarrin@udla.edu.ec}}

\affil[1]{\scriptsize Escuela de Ciencias Físicas y Matemáticas, Universidad de Las Américas, Vía a Nayón, C.P.170124, Quito, Ecuador.} 
\date{\today}
\begin{document}
	\maketitle	 
	\begin{abstract}
Uniqueness of the trivial solution (the zero solution) for the steady-state Navier-Stokes equations is an interesting problem who has known several recent contributions.  These results are  also known as the Liouville type problem  for the steady-state Navier-Stokes equations. In the setting of the $L^p-$  spaces, when $3\leq p \leq 9/2$ it is known that the trivial solution of these equations is the unique one. In this note, we extend this previous result to other values of the  parameter $p$. More precisely, we prove that the velocity field must be zero provided that it belongs to the $L^p -$ space with $3/2<p<3$. Moreover, for the large interval of values $9/2<p<+\infty$, we also obtain a partial result on the vanishing of the velocity  under an additional  hypothesis in terms of the  Sobolev space of negative order $\dot{H}^{-1}$. This last result  has an interesting corollary when studying   the Liouville problem in the natural energy space of these solutions $\dot{H}^{1}$.  \\[3mm]
	\end{abstract}

\section{Introduction} 
This short note deals with the  homogeneous and incompressible  steady-state (time-independent) Navier-Stokes  equations in the whole space $\Rt:$ 
\begin{equation*}
(NS) \quad -\Delta \U+(\U \cdot \vec{\nabla})\U +\vec{\nabla}P = 0, \qquad div(\U)=0.
\end{equation*}    
Here $\U=(U_1, U_2, U_3): \Rt\longrightarrow \Rt$ is the velocity of the fluid while   $P:\Rt\longrightarrow \mathbb{R}$ is the pressure. The  equation $\ds{div(\U)=0}$ represents the fluid incompressibility.  We recall that $(\U,P)$ is a smooth solution for the (NS) equations if $\U \in \mathcal{C}^{2}(\Rt)$, $P \in \mathcal{C}^{1}(\Rt)$ and if the couple $(\U,P)$ verify these equations in the classical sense. 

\medskip

In the  (NS) equations,  we may observe that $\U=0$ and $P=0$ is always a smooth solution, also known as the \emph{trivial solution}, and then it is quite natural to ask if the trivial solution is the unique one. In the general setting of the space $\mathcal{C}^{2}(\Rt) \times \mathcal{C}^{1}(\Rt)$, the answer to this question is negative and we are able to give a simple counterexample. Let $f \in \mathcal{C}^{3}(\Rt)$ be a scalar field such that $\Delta f =0$. Then, by setting the velocity $\U= \vec{\nabla} f$ and the pressure $P= - \frac{1}{2} \vert \vec{\nabla} f \vert^2$, and by using some well-known rules of the vector calculus, we have that $(\U,P)$ is also a smooth solution of the (NS) equations. See the Appendix \ref{AppendixA}  for all the details.  We thus look for some \emph{additional hypothesis} on  smooth solutions of the (NS) equations to insure the uniqueness of the trivial solution. This type of problem is also  known as a the \emph{Liouville problem} for the  (NS) equations. We emphasize this problem  has attired a lot of attention in the community of researchers. See, for instance, \cite{ChaeWolf,ChaeWeng,OscarPaper,Ser2016,Seregin3,SereWang} and the references therein. 

\medskip

In the example above, we remark that as the scalar field $f$ is a harmonic function then it is a polynomial and thus $\U = \vec{\nabla} f$ also has a polynomial growth at infinity. This fact strongly suggests that, in order to study the Liouville problem for the (NS) equations, we must seek for \emph{decaying properties} at infinity on the velocity $\U$. Precisely, it  was pointed out in the celebrated work of  G. Galdi who showed in \cite{Galdi} (Chapter X,  Remark X.9.4 and Theorem X.9.5, p. 729) that if $\U$ is a smooth solution of the (NS) equations, and moreover, if $\U \in L^{9/2}(\Rt)$ then we necessarily have $\U=0$.   This result follows from the following \emph{Cacciopoli} type estimate: for $R>1$, we denote    $B_R=\{ x \in \Rt: \vert x \vert <R \}$ and  $C(R/2,R)=\{ x \in \Rt: R/2 < \vert x \vert <R \}$, and for a constant $C>0$ which does not depend on $R$  we have:  
\begin{equation*}
\int_{B_R} \vert \vec{\nabla} \otimes \U (x) \vert^2 dx \leq \,  C \Vert \U \Vert^{3}_{L^{9/2}(C(R/2,R))} + C  R^{-1/3} \Vert \U \Vert^{2}_{L^{9/2}(C(R/2,R))} + C \Vert \U \Vert_{L^{9/2}(C(R/2,R))} \Vert P \Vert_{L^{9/4}(C(R/2,R))}.
\end{equation*}
This estimate  yields the identity  $\U=0$, provided that $\U \in L^{9/2}(\Rt)$  and   $P \in L^{9/4}(\Rt)$. Moreover, it is worth mentioning  the value of the integration parameter $9/2$ naturally appears by the well-known scaling properties of (NS) equations.  

\medskip

Galdi's result was recently extended in \cite{Liouville1} for other values of the parameter $p$ in the Lebesgue spaces. More precisely, in Theorem $1$ of \cite{Liouville1}, D. Chamorro, P.G. Lemarié-Rieusset and the author of this note proved that if $\U$ is a smooth solution of the (NS) equations  such that  $\U \in L^{p}(\Rt)$, with $3 \leq p < 9/2$, then we have $\U=0$ and $P=0$. For this, we performed   the next  \emph{Cacciopoli} type  estimate: 
\begin{equation*}
\int_{B_R} \vert \vec{\nabla} \otimes \U (x) \vert^2 dx \leq C\,   R^{1-6/p} \Vert \U \Vert_{L^p(B_R)}+ C_p\,  R^{2-9/p} \Vert \U \Vert^{3}_{L^p(C(R/2,R))} + C_p \, \Vert \U \Vert_{L^p(C(R/2,R))}\, \Vert P \Vert_{L^{p/2}(C(R/2,R))}. 
\end{equation*}
Here, the constant $C_p>0$ is the norm of a certain test function in the space $L^{\frac{p}{p-3}}(\Rt)$ and we thus need condition $ 3 \leq p$. On the other hand, in order to get a uniform control on $R>1$, due to the expression $R^{2-9/p}$ in the second term on  the right-hand side, we also need the condition $p\leq 9/2$. Consequently, the Liouville problem is solved in the $L^p-$ space for $3 \leq p \leq 9/2$.  

\medskip

The aim of this short note is to continue  with the  study  of the Liouville problem for the (NS) equations in the setting of the Lebesgue spaces. Specifically, we study this problem   for the  values of the parameter $p$ outside the interval $3 \leq p \leq 9/2$.  For the values of the parameter $p$ lower than $3$, our first result states as follows. 
\begin{TheoremeP}\label{Th1} Let $(\U, P)$ be a smooth solution of the (NS) equations. If $\U \in L^{p}(\Rt)$, with $\frac{3}{2}<p<3$, then we have  $\U=0$ and $P=0$.
\end{TheoremeP}	

The proof is based on a different and more technical \emph{Cacciopoli} type estimate, which is  stated in  Proposition \ref{lemaTehc-Ser}.  This  estimate was  originally established in \cite{Seregin3}. It is worth mentioning our computations are not longer valid when  $1\leq p \leq 3/3$ and, to the best of our knowledge, the Liouville problem for the (NS) equations is still an open question  these values of the parameter $p$.

\medskip 

It is also worth emphasizing some recent preprints study this problem when $\U \in L^p(\Rt)$, whit $1 \leq p  < 3$, provided that $\vert \U(x) \vert \to 0$ when $\vert x \vert \to +\infty$. However, as $\U$ is assumed a smooth solution, the last hypothesis implies that $\U \in L^{\infty}(\Rt)$. Then, by the interpolations inequalities we have $\U \in L^p \cap L^{\infty}(\Rt)$, hence $\U \in L^{9/2}(\Rt)$. Consequently, the additional vanishing  assumption at infinity makes the problem trivial. In the sense, one the main interests on the result above is the study of the Liouville problem in the interval $3/2<p<3$ without any additional assumption. 

\medskip

When $9/2 <p <6$, the Liouville  problem in the $L^p -$ spaces was study  studied by using  some  supplementary hypothesis.  In \cite{Liouville1}, this problem is solved in the space $L^{p} \cap B^{3/p-3/2}_{\infty, \infty} (\Rt)$,  where $B^{3/p-3/2 }_{\infty,\infty} (\Rt)$ is a homogeneous Besov space (\ref{Besov-neg}). Moreover,  for the value $p=6$, an interesting result of G. Seregin given in \cite{Ser2016} shows that this problem is  solved in the space $L^{6}\cap BMO^{-1}(\Rt)$. On the other hand, to the best of our knowledge, there are not previous results for the values $6<p$. In our next result, we  study the Liouville problem in the large interval $9/2<p<+\infty$. As the previous results, we shall need here an additional assumption on the velocity field.

\begin{TheoremeP}\label{Th2} Let $(\U, P)$ be a smooth solution of the (NS) equations.  If  $\U \in L^{p}(\Rt)\cap \dot{H}^{-1}(\Rt)$,  for $9/2 < p <+\infty$,  then we have $\U=0$ and $P=0$.
\end{TheoremeP}	

One the main interests of this result is the use of the space $\dot{H}^{-1}(\Rt)$, which is the dual space of  functions having a finite Dirichlet integral; and which provides us a different condition to solve the Liouville type problem. On the other hand, this result suggests that the value $p=9/2$, found by G. Galdi in \cite{Galdi}, seems to be the upper limit  to solve the Liouville problem without any additional assumptions.  This fact also suggests to look for non trivial solutions for the (NS) equations in the space $L^p(\Rt)$ with $9/2<p$. However, we think that this is still a very challenging open question.  

\medskip

In this result, the  value $p=6$ is of particular interest. Let us recall that the problem Liouville for the (NS) equations  was originally stated in the   homogeneous Sobolev space $\dot{H}^{1}(\Rt)$ (see for instance   the Chapter $X$ in \cite{Galdi} and the Chapter $4$ in \cite{OscarTesis}). This is the natural \emph{energy space} for this equations, and moreover, not rigorous computations strongly suggest that all the solutions $\U \in \dot{H}^{1}(\Rt)$ must be identical to zero. 
 However,  the information $\U \in \dot{H}^{1}(\Rt)$ is not enough to justify all the computations, in particular those estimates involving the nonlinear term $\ds{(\U \cdot \vec{\nabla}) \U}$. Consequently,  the Liouville  problem for the (NS) equations in the space $\dot{H}^{1}(\Rt)$ is still an outstanding  open question far from obvious. 

\medskip

By the Sobolev embeddings we have $\dot{H}^{1}(\Rt)\subset L^6(\Rt)$, and then, we can use the Theorem \ref{Th2} to prove the following corollary. This is  a partial result on the Liouville problem in the space $\dot{H}^{1}(\Rt)$.  

\begin{CorollaireP}\label{Col} Let be $\U \in \dot{H}^{1}(\Rt)$ be a  weak solution of the (NS) equations.   If $\ds{\widehat{\U} \in L^{r}_{loc}(\Rt)}$, with $r>6$, then we have $\U =0$. 
\end{CorollaireP}	

This corollary does not assume any  smoothness of the velocity $\U$  and we understand here $\U \in \dot{H}^{1}(\Rt)$ as a weak solution, \emph{i.e.}, $\U$ verifies the (NS) equations in the distributional sense. On the other hand, we recall  that the space  $\dot{H}^{1}(\Rt)$ is defined as the space of temperate distributions $g$ such that $\widehat{g} \in L^{1}_{loc}(\Rt)$ and $\ds{\int_{\Rt} \vert \xi \vert^2 \vert \widehat{g}(\xi) \vert^2 \, d \xi <+\infty}$.  We thus  observe that a stronger locally-integrability condition on  $\widehat{\U}$ yields the desired identity $\U=0$.



\section{The key tools}\label{Sec2:Previous-Results}  
\subsection{Homogeneous Besov spaces.}
The first key tool deals with the  homogeneous Besov spaces. Let $0<s<1$, and $1\leq p,q\leq +\infty$. The homogeneous Besov space of positive order: $\dot{B}^{s}_{p, q}(\Rt)$, is defined as the set of $f \in \mathcal{S}'(\Rt)$ such that  
\begin{equation}\label{besov-pos1}
\Vert f \Vert_{\dot{B}^{s}_{p, q}}=  \left( \int_{\Rt} \frac{\Vert f(\cdot+x)-f(\cdot) \Vert_{L^p}}{\vert x \vert^{3+sq}}\, dx \right)^{1/q} <+\infty, \quad \text{with}, \quad 1\leq p,q<+\infty,
\end{equation}
and 
\begin{equation}\label{besov-pos}
\Vert f \Vert_{\dot{B}^{s}_{\infty, \infty}}= \sup_{x \in \Rt} \frac{\Vert f(\cdot+x)-f(\cdot)\Vert_{L^{\infty}}}{\vert x \vert^{s}}<+\infty.
\end{equation}
Moreover, the Besov space of negative order $\dot{B}^{-s}_{\infty, \infty}(\Rt)$ can be characterized by means of the heat kernel $h_t$  as the set of $f \in \mathcal{S}'(\Rt)$ such that 
\begin{equation}\label{Besov-neg}
\Vert f \Vert_{\dot{B}^{-s}_{\infty, \infty}} = \sup_{t>0} \, t^{\frac{s}{2}} \Vert h_t \ast f \Vert_{L^{\infty}}<+\infty.
\end{equation}
For more details on the Besov spaces  and their application to the theoretical study of the Navier-Stokes equations (stationary or time-dependent)  see the Chapter $8$ in the book \cite{PGLR1}.

\subsection{A \emph{Cacciopoli} type estimate} 

The second key tool is the following \emph{Cacciopoli} type estimate.  This inequality is inspired by the  estimate $(2.2)$, page $9$ in \cite{Seregin3} and its proof essentially follows similar ideas. However, for the sake of completeness, we include the main steps of the proof. 

\begin{Proposition}\label{lemaTehc-Ser} Let $(\U, P)$ be a smooth solution of the (NS) equation. Let $\V_1$ and $\V_2$ be smooth vector fields such that $\U = \vec{\nabla} \wedge \V_1$ and $\U= \vec{\nabla} \wedge \V_2$. Then, for all $3<q<+\infty $, there exists a constant $C_q>0$ such that such that for all $R>1$ we have: 
		\begin{equation*}
		\int_{B_{R/2}} \vert \vec{\nabla} \otimes \U(x) \vert^2 dx \leq  \frac{C_q}{R} \left( \frac{1}{R^3} \int_{\mathcal{C}(R/2,R)} \vert \V_1(x) \vert^2 dx \right)  \left( 1+ \left(\frac{1}{R^3} \int_{B_R} \vert \V_2(x) \vert^{q}dx \right)^{\frac{4}{q-3}} \right).
		\end{equation*}
	\end{Proposition}	  
\pv We introduce the following cut-off functions, which were considered in \cite{OscarPaper}.  For  $R>1$ fixed,  we define first the function $\varphi_R\in \mathcal{C}^{\infty}_{0}(B_R)$ such that: for $R/2 < \rho < r< R$ it verifies: $\varphi_R(x)=1$ when $\vert x \vert < \rho$, $\varphi_R(x)=0$ when $\vert x \vert >r$  and, for all multi-indice $\vert \alpha \vert \leq 4$, $\Vert \partial^{\alpha} \varphi_R \Vert_{L^\infty} \leq \frac{c}{(r-\rho)^{\vert \alpha \vert}}$. Next, we define the function $\W_R$ as the solution of the problem:
\begin{equation*}
div(\W_R)=\vec{\nabla}\varphi_R\cdot \U, \quad \text{in}\,\, B_r, \quad \text{and}\quad \W_R=0 \,\, \text{on}\,\,  \partial B_r. 
\end{equation*} 
In Lemma  $III. 3.1$,  page 162 of  \cite{Galdi}, it is proven that for $1<q<+\infty$ there exists a solution $\W_R\in W^{1,q}(B_R)$, which verifies  $supp\,(\W_R)\subset \mathcal{C}(R/2, R)$  and $\Vert \vec{\nabla}\otimes \W_R\Vert_{L^q(\mathcal{C}(R/2, R))}\leq c \Vert   \vec{\nabla}\varphi_R\cdot \U \Vert_{L^q(\mathcal{C}(R/2, R))}$.

\medskip 

We consider now  the function $\varphi_R \U-\W_R$, which has a support in the ball $B_r$,  and we write  
\begin{equation}\label{eq01}
\int_{B_r} \left( -\Delta \U +(\U \cdot \vec{\nabla})\U +\vec{\nabla}P\right)\cdot \left( \varphi_R \U-\W_R \right)dx=0.
\end{equation}
From this identity, by performing some integration by parts,  and moreover, for a function $\V_2$ such that $\U = \vec{\nabla} \wedge \V$ on $B_R$, for $1=1/p+1/q$ we have the following estimate (see the page $6$ of \cite{Seregin3} or the appendix of \cite{Liouville1}): 
\begin{equation*}
\begin{split}
\int_{B_\rho} \vert \vec{\nabla} \otimes \U \vert^2 dx  \leq & \, C \left( \int_{B_r} \vert \vec{\nabla} \otimes \U \vert^2 dx\right)^{1/2} \left( \int_{B_r} \vert \vec{\nabla} \varphi_R \vert^2 \vert \U \vert^2 dx\right)^{1/2}\\
&\, + C \left(\int_{B_r} \vert \vec{\nabla} \otimes \U \vert^2 dx\right)^{1/2} \left( \int_{B_r}\vert \vec{\nabla} \varphi_R \vert \vert \U \vert \vert \V_2 \vert dx \right)^{1/2}\\
\leq & \, C \left( \int_{B_r} \vert \vec{\nabla} \otimes \U \vert^2 dx\right)^{1/2} \left( \int_{B_r} \vert \vec{\nabla} \varphi_R \vert^2 \vert \U \vert^2 dx\right)^{1/2}\\
&\, + C \left(\int_{B_r} \vert \vec{\nabla} \otimes \U \vert^2 dx\right)^{1/2} \left( \int_{B_r}\vert \vec{\nabla} \varphi_R \vert^p \vert \U \vert^p dx \right)^{1/p} \left( \int_{B_r} \vert \V_2 \vert^q dx \right)^{1/q}.
\end{split}
\end{equation*} 
Then, we  apply the discrete Young inequalities (with $1=1/2+1/2$) in both terms on the right-hand side to get:
\begin{equation}\label{estim-01}
\int_{B_\rho} \vert \vec{\nabla} \otimes \U \vert^2 dx  \leq  \frac{1}{8} \int_{B_r} \vert \vec{\nabla} \otimes \U \vert^2 dx + C \int_{B_r} \vert \vec{\nabla} \varphi_R \vert^2 \vert \U \vert^2 dx + C \left( \int_{B_r}\vert \vec{\nabla} \varphi_R \cdot  \U \vert^p dx \right)^{2/p} \left( \int_{B_r} \vert \V_2 \vert^q dx \right)^{2/q}.
\end{equation}
We must study the third term on the right-hand side. We set $2<p<6$. For $\theta= \frac{3(p-2)}{2p}$ (which verifies $0<\theta<1$) by the interpolation inequalities (with $2/p=\theta/3 + (1-\theta)/1$) and the Sobolev embedding $ \left(\int_{B_r} \vert f \vert^6 dx \right)^{1/6} \leq c \left( \int_{B_r} \vert \vec{\nabla} f \vert^2 dx\right)^{1/2}$,  we write: 
\begin{equation*}
\begin{split}
&\left( \int_{B_r}\vert \vec{\nabla} \varphi_R \cdot \U \vert^p dx \right)^{2/p}  = \,  \left( \int_{B_r} \left(\vert \vec{\nabla} \varphi_R  \cdot \U \vert^2 \right)^{p/2} dx \right)^{2/p} \\
  \leq & \,  C \left( \int_{B_r}\vert \vec{\nabla} \varphi_R \vert^2\, \vert \U \vert^2 \, dx \right)^{1-\theta}\left( \int_{B_r} \left( \vert \vec{\nabla} \varphi_R \cdot  \U \vert^2 \right)^3 dx \right)^{\theta} \\
 \leq & \, C \left( \int_{B_r}\vert \vec{\nabla} \varphi_R \vert^2\, \vert \U \vert^2 \, dx \right)^{1-\theta} \left( \int_{B_r} \left\vert \vec{\nabla} \left(  \vec{\nabla} \varphi_R \cdot \U   \right)  \right\vert^2 dx \right)^{\theta} \\
 \leq & \, C  \left( \int_{B_r}\vert \vec{\nabla} \varphi_R \vert^2\, \vert \U \vert^2 \, dx \right)^{1-\theta} \left[ \left( \int_{B_r} \left\vert  \vec{\nabla} ( \vec{\nabla} \  \varphi_R  )    \right\vert^2\, \vert \U \vert^2 dx \right)^{\theta} + \left( \int_{B_r} \vert \vec{\nabla} \varphi_R \vert^2 \vert \vec{\nabla} \otimes \U \vert^2 dx\right)^{\theta} \right]\\
 \leq & \, C \left( \int_{C(\rho, r)}\vert \vec{\nabla} \varphi_R \vert^2\, \vert \U \vert^2 \, dx \right)^{1-\theta} \left[ \left( \int_{C(\rho, r)} \left\vert  \vec{\nabla} ( \vec{\nabla} \  \varphi_R  )    \right\vert^2\, \vert \U \vert^2 dx \right)^{\theta} + \left( \int_{B_r} \vert \vec{\nabla} \varphi_R \vert^2 \vert \vec{\nabla} \otimes \U \vert^2 dx\right)^{\theta} \right]
\end{split}
\end{equation*}
With this estimate, we get back to (\ref{estim-01}). We applying again the discrete Young inequalities (with $1=\theta + (1-\theta)$) and moreover, by the estimate $\Vert \Delta \varphi_R \Vert_{L^\infty} \leq \frac{c}{(r-\rho)^2}$,   we write: 
\begin{equation}\label{estim-02}
\begin{split}
\int_{B_\rho} \vert \vec{\nabla} \otimes \U \vert^2 dx  \leq & \, \frac{1}{4} \int_{B_r} \vert \vec{\nabla} \otimes \U \vert^2 \, dx  + C \underbrace{\int_{C(\rho,r)}\vert \vec{\nabla} \varphi_R \vert^2\, \vert \U \vert^2 \, dx}_{(\mathcal{A})} \\
 & \, + C \left(  \int_{C(\rho,r)} \vert \vec{\nabla} \varphi_R \vert^2 \vert \U \vert^2 dx\right)^{1-\theta} \underbrace{ \left( \int_{C(\rho,r)} \left\vert\vec{\nabla} \left(\vec{\nabla} \varphi_R\right) \right\vert^2 \vert \U \vert^2 dx \right)^{\theta}}_{(\mathcal{B})} \left( \int_{B_r} \vert \V_2 \vert^q dx\right)^{2/q}\\
 &\, + C \left(\frac{1}{(r-\rho)^2}\right)^{\frac{\theta}{1-\theta}}\, \left( \int_{C(\rho,r)} \vert \vec{\nabla} \varphi_R \vert^2\, \vert \U \vert^2 dx\right)\left( \int_{B_r} \vert \V_2 \vert^{q} dx\right)^{\frac{2}{q(1-\theta)}}.
\end{split}
\end{equation}
We must estimate now the terms $(\mathcal{A})$ and $(\mathcal{B})$. By Lemma $2.1$ of \cite{Seregin3}, for a smooth function $\V_1$ such that $\U = \vec{\nabla} \wedge \V_1$ we have:
\[ (\mathcal{A}) \leq \frac{C}{(r-\rho)^2} \underbrace{ \left[  \left(  \int_{B_r} \vert \vec{\nabla} \otimes \U \vert^2 dx \right)^{1/2} \left( \int_{C(\rho,r)} \vert \V_1 \vert^2 dx\right)^{1/2} + \frac{1}{(r-\rho)^2} \int_{C(\rho,r)} \vert \V_1 \vert^2 dx \right]}_{(\mathcal{C})}, \]
\[ (\mathcal{B}) \leq \frac{C}{(r-\rho)^4} \underbrace{\left[  \left(  \int_{B_r} \vert \vec{\nabla} \otimes \U \vert^2 dx \right)^{1/2} \left( \int_{C(\rho,r)} \vert \V_1 \vert^2 dx\right)^{1/2} + \frac{1}{(r-\rho)^2} \int_{C(\rho,r)} \vert \V_1 \vert^2 dx \right]}_{(\mathcal{C})}. \]
We get back to the estimate (\ref{estim-02}) to write:
\begin{equation*}
\begin{split}
 \int_{B_\rho} \vert \vec{\nabla} \otimes \U \vert^2 dx  \leq  & \, \frac{1}{4} \int_{B_r} \vert \vec{\nabla} \otimes \U \vert^2 dx + \frac{C}{(r-\rho)^2} (\mathcal{C}) \\
& \, + \left(  \frac{C}{(r-\rho)^2}\, (\mathcal{C})  \right)^{1-\theta} \left(\frac{C}{(r-\rho)^4}\, (\mathcal{C}) \right)^{\theta}\, \left( \int_{B_r} \vert \V_2 \vert^{q} dx\right)^{2/q}  \\
&\, + \left( \frac{C}{(r-\rho)^2} \right)^{\frac{\theta}{1-\theta}}\, \frac{C}{(r-\rho)^2}(\mathcal{C}) \left( \int_{B_r} \vert \V \vert^q\, dx\right)^{\frac{2}{q(1-\theta)}}\\
= &\,  \frac{1}{4} \int_{B_r} \vert \vec{\nabla} \otimes \U \vert^2 dx + \frac{C}{(r-\rho)^2} (\mathcal{C}) + \frac{C}{(r-\rho)^{2 \theta}}\, \frac{1}{(r-\rho)^2} (\mathcal{C}) \left( \int_{B_r} \vert \V_2 \vert^q \, dx\right)^{2/q} \\
&\, + \left( \frac{C}{(r-\rho)^2} \right)^{\frac{\theta}{1-\theta} +1} \, (\mathcal{C})\, \left( \int_{B_r} \vert \V_2 \vert^q\, dx\right)^{\frac{2}{q(1-\theta)}}.
\end{split}
\end{equation*}
In the last term above, we remark that we have $\frac{\theta}{1-\theta} + 1 = \frac{1}{1-\theta}$. Moreover,  by definition of the parameter $\theta= \frac{3(p-2)}{2p}$  and by the relation $1/2 = 1/p +1/q$ we have $q\theta =3$. Thus, we can write:
\begin{equation*}
\begin{split}
&\int_{B_\rho} \vert \vec{\nabla} \otimes \U \vert^2 dx \\
 \leq & \,   \frac{1}{4} \int_{B_r} \vert \vec{\nabla} \otimes \U \vert^2 dx + \frac{C}{(r-\rho)^2} (\mathcal{C})\, \left[ 1 + \left( \frac{1}{(r-\rho)^{q \theta}} \int_{B_r} \vert \V_2 \vert^q\, dx\right)^{2/q} + \left( \frac{1}{(r-\rho)^{q \theta}} \int_{B_r} \vert \V_2 \vert^q\, dx\right)^{\frac{2}{q(1-\theta)}}\right]\\
=&\,  \frac{1}{4} \int_{B_r} \vert \vec{\nabla} \otimes \U \vert^2 dx + \frac{C}{(r-\rho)^2} (\mathcal{C})\, \underbrace{\left[ 1 + \left( \frac{1}{(r-\rho)^{3}} \int_{B_r} \vert \V_2 \vert^q\, dx\right)^{2/q}  + \left( \frac{1}{(r-\rho)^{3}} \int_{B_r} \vert \V_2 \vert^q\, dx\right)^{\frac{2}{q(1-\theta)}}\right]}_{(\mathcal{D})}.
\end{split}
\end{equation*}
Now, we write down the whole term $(\mathcal{C})$, and we have  
\begin{equation*} 
\begin{split}
&\int_{B_\rho} \vert \vec{\nabla} \otimes \U \vert^2 dx \\
 \leq &\, \frac{1}{4} \int_{B_r} \vert \vec{\nabla} \otimes \U \vert^2 dx + \frac{C}{(r-\rho)^2} \left[ \left(  \int_{B_r} \vert \vec{\nabla} \otimes \U \vert^2 dx \right)^{\frac{1}{2}} \left( \int_{C(\rho,r)} \vert \V_1 \vert^2 dx\right)^{\frac{1}{2}} + \frac{1}{(r-\rho)^2} \int_{C(\rho,r)} \vert \V_1 \vert^2 dx  \right](\mathcal{D}) \\
\leq &\, \frac{1}{4} \int_{B_r} \vert \vec{\nabla} \otimes \U \vert^2 dx+ \frac{C}{(r-\rho)^2} \left[ \frac{(r-\rho)^2}{4C} \int_{B_r} \vert \vec{\nabla} \otimes \U \vert^2 dx + \frac{C}{(r-\rho)^2} \int_{C(\rho,r)} \vert \V_1 \vert^2 \, dx \, (\mathcal{D})^2 \right]\\
\leq &\, \frac{1}{2} \int_{B_r} \vert \vec{\nabla}\otimes \U \vert^2 dx + \frac{C}{(r-\rho)^4} \int_{C(\rho,r)} \vert \V_1 \vert^2 dx \, (\mathcal{D})^2
\end{split}
\end{equation*} 
Now, we write down the whole term $(\mathcal{D})$ to write:
\begin{equation*}
\begin{split}
\int_{B_\rho} \vert \vec{\nabla} \otimes \U \vert^2 dx 
\leq &\,\frac{1}{2} \int_{B_r} \vert \vec{\nabla} \otimes \U \vert^2 dx + \frac{C}{(r-\rho)} \left[ \frac{1}{(r-\rho)^3}  \int_{C(\rho,r)}\vert \V_1 \vert^2 dx \right] \\
&\, \times   \left[1 + \left( \frac{1}{(r-\rho)^3} \int_{B_r} \vert \V_2 \vert^q\, dx\right)^{4/q} + \left( \frac{1}{(r-\rho)^3} \int_{B_r} \vert \V_2 \vert^q dx\right)^{\frac{4}{q(1-\theta)}} \right].  
\end{split}
\end{equation*} 
We remark that as $0<1-\theta <1$ and as $q\theta=3$, we have  $\frac{4}{q} \leq \frac{4}{q(1-\theta)}= \frac{4}{q-3}$. Then, from the last estimate we obtain: 
\begin{equation*}
\begin{split}
\int_{B_\rho} \vert \vec{\nabla} \otimes \U \vert^2 dx \leq\, \frac{1}{2} \int_{B_r} \vert \vec{\nabla} \otimes \U \vert^2 dx+\frac{C}{(r-\rho)} \left[ \frac{1}{(r-\rho)^3}  \int_{C(\rho,r)}\vert \V_1 \vert^2 dx \right] \, \left[ 1 + \left( \frac{1}{(r-\rho)^3} \int_{B_r} \V_2 \vert^q \, dx\right)^{\frac{4}{q-3}}\right].
\end{split}
\end{equation*}
Finally, we use a well-known iterative argument, see for instance  \cite{Liouville1} and \cite{Giaquinta}, to obtain the Cacciopoli type estimate stated in Proposition \ref{lemaTehc-Ser}, which is now proven.   \finpv 
\section{Proofs of the results}\label{Sec3:Proofs}
\subsection{Proof of Theorem \ref{Th1}}
Let $(\U, P)$ be a smooth solution of the (NS) equations. We assume that $\U \in L^{p}(\Rt)$ with $3/2<p<3$. In the framework of  Proposition \ref{lemaTehc-Ser}, we will set the vector fields $\V_1$ and $\V_2$ as follows:  we define first the vector field $\V$ by means of the velocity $\U$ as $\ds{\V= \frac{1}{-\Delta} (\vec{\nabla} \wedge \U),}$  where we have $\U= \vec{\nabla} \wedge \V$. Indeed, as  we have $div(\U)=0$, then we can write 
	$$ \vec{\nabla}\wedge \V=\vec{\nabla}\wedge \left( \vec{\nabla}\wedge \left( \frac{1}{-\Delta} \U \right) \right)=\vec{\nabla} \left( div \left( \frac{1}{-\Delta} \U  \right)\right) - \Delta \left( \frac{1}{-\Delta} \U  \right) =\U.$$
	Now, for all $x \in \Rt$ we set the vector fields $\V_1(x)= \V(x)$ and $\V_2(x)= \V(x)$. Then,  for all $3<q<+\infty$ and for all $R>1$, by Proposition \ref{lemaTehc-Ser} we have the estimate: 
	\begin{equation}\label{eq03}
	\int_{B_{R/2}} \vert \vec{\nabla} \otimes \U(x) \vert^2 dx \leq  \frac{C_q}{R} \left( \frac{1}{R^3} \int_{\mathcal{C}(R/2,R)} \vert \V(x) \vert^2 dx \right)  \left( 1+ \left(\frac{1}{R^3} \int_{B_R} \vert \V(x) \vert^{q}dx \right)^{\frac{4}{q-3}} \right).
	\end{equation}
	We must study now the term on the right-hand side. We  recall that we have $\ds{\V= \frac{1}{-\Delta} (\vec{\nabla} \wedge \U),}$ and moreover, we have $\U \in L^{p}(\Rt)$ with $3/2<p<3$. Then,  we  get  $\V \in L^{3p/(3-p)}(\Rt)$ with $3p/(3-p)>3$. Indeed,  we write $\ds{\V= \frac{1}{\sqrt{-\Delta}} (\frac{1}{\sqrt{-\Delta}} (\vec{\nabla} \wedge \U))}$, and then, by the properties of the Riesz potential $\frac{1}{\sqrt{-\Delta}}$, as well  by the properties of the Riesz transforms $\frac{\partial_i}{\sqrt{-\Delta}}$, we can write: 
	$$  \Vert \V \Vert_{L^{3p/(3-p)}} \leq c \left\Vert  \frac{1}{\sqrt{-\Delta}} (\frac{1}{\sqrt{-\Delta}} (\vec{\nabla} \wedge \U)) \right\Vert_{L^{3p/(3-p)}} \leq c \left\Vert  \frac{1}{\sqrt{-\Delta}} (\vec{\nabla} \wedge \U) \right\Vert_{L^p}  \leq c \Vert \U \Vert_{L^p}. $$
	Moreover, as $3/2<p<3$ we have $3p/(3-p)>3$. Thus, we set the parameter  $q=3p/(3-p)$ and we have $\Vert \V \Vert_{L^q}<+\infty$.  
	
	\medskip 
	
	We get back to the estimate  (\ref{eq03}) to write:  
\begin{equation*}
\begin{split}
	\int_{B_{R/2}} \vert \vec{\nabla} \otimes \U(x) \vert^2 dx \leq & \,  \frac{C_q}{R^4} \left(  \int_{\mathcal{C}(R/2,R)} \vert \V(x) \vert^2 dx \right)  \left( 1+ \left(\frac{1}{R^3} \int_{B_R} \vert \V(x) \vert^{q}dx \right)^{\frac{4}{q-3}} \right) \\
	\leq & \,  \frac{C_q}{R^4} \,  R^{6(1/2-1/q)} \left( \int_{\mathcal{C}(R/2,R)} \vert \V(x)\vert^q dx \right)^{2/q} \left( 1+ \left(\frac{1}{R^3} \int_{B_R} \vert \V(x) \vert^{q}dx \right)^{\frac{4}{q-3}} \right)\\
	\leq & \, C_q R^{-1-6/q} \left( \int_{\mathcal{C}(R/2,R)} \vert \V(x)\vert^q dx \right)^{2/q} \left( 1+ \left(\frac{1}{R^3} \int_{B_R} \vert \V(x) \vert^{q}dx \right)^{\frac{4}{q-3}} \right).
\end{split}
\end{equation*}	
Now, we let $R\to +\infty$ and we obtain $\ds{\int_{\Rt} \vert \vec{\nabla} \otimes \U(x) \vert^2 dx =0}$.  But, by the Sobolev embeddings we write  $\Vert \U \Vert_{L^6}\leq  c \Vert \vec{\nabla} \otimes \U \Vert_{L^2}$ and then we have the identity $\U=0$. Finally, by splitting  the pressure $P$ as the well-known expression $\ds{P =\sum_{1 \leq i,j \leq 3} \mathcal{R}_{i}\mathcal{R}_{j}(U_i U_j)}$, we  conclude that $P=0$. Theorem \ref{Th1} is now proven.  \finpv
\subsection{Proof of Theorem \ref{Th2}} 
	
 We assume now that the smooth solution $(\U,P)$ verifies  $\U \in L^{p} \cap \dot{H}^{-1}(\Rt)$ with $9/2 \leq p < +\infty$. As before, we define  the vector field $\ds{\V= \frac{1}{-\Delta} (\vec{\nabla} \wedge \U)}$, and we set now the vector fields $\V_1$ and $\V_2$ as follows: 
	\begin{equation}\label{def-V1-V2}
	\V_1(x)= \V(x) \quad \text{and} \quad \V_2(x)= \V(x)-\V(0). 
	\end{equation}
	We remark that we have $\U = \vec{\nabla} \wedge \V_1$ and $\U= \vec{\nabla} \wedge \V_2$, and then, for $q=p$ and  for all $R>1$, by Proposition  \ref{lemaTehc-Ser}  we can write: 
	\begin{equation}\label{Caccioppoli2}
	\begin{split}
	\int_{B_{R/2}} \vert \vec{\nabla} \otimes \U(x) \vert^2 dx 
	\leq & \, \frac{C_p}{R} \left(\frac{1}{R^3}\int_{\mathcal{C}(R/2,R) } \vert \V_1(x) \vert^2 dx \right) \left( 1+ \left(\frac{1}{R^3} \int_{B_R} \vert \V_2(x) \vert^{p}dx 
	\right)^{\frac{4}{p-3}} \right) \\
	\leq & \, frac{C_p}{R^4}   \left(\int_{\mathcal{C}(R/2,R)} \vert \V_1(x) \vert^2 dx \right) \left( 1+ \left(\frac{1}{R^3} \int_{B_R} \vert \V_2(x) \vert^{p}dx 
	\right)^{\frac{4}{p-3}} \right) \\
	\leq & \, C_p  \left( \int_{\mathcal{C}(R/2,R)} \vert \V_1(x) \vert^2 dx \right)\,  \left(\frac{1}{R^4}+ \underbrace{\frac{1}{R^4} \left(\frac{1}{R^3} \int_{B_R} \vert \V_2(x) \vert^{p}dx 
		\right)^{\frac{4}{p-3}}}_{(I(R))}  \right).
\end{split}
	\end{equation}
	
We must estimate the term $I(R)$. 	As we have $\U \in L^p(\Rt)$, and moreover, by  the continuous embedding $L^p(\Rt) \subset \dot{B}^{-3/p}_{\infty, \infty}(\Rt)$, we obtain $\U \in \dot{B}^{-3/p}_{\infty, \infty}(\Rt)$. Then,  as  $\ds{\V= \frac{1}{-\Delta} (\vec{\nabla} \wedge \U)}$  then we get $\V \in \dot{B}^{1-3/p}_{\infty, \infty}(\Rt)$. Moreover, as we have $9/2 \leq p < +\infty$, then we get  $\frac{1}{3}\leq 1-\frac{3}{p}<1$, where  the Besov space $\dot{B}^{1-3/p}_{\infty, \infty}(\Rt) $ is defined in the formula (\ref{besov-pos}). By   this formula, for all $R>1$,  we can write $$\ds{\sup_{x \in B_R} \frac{\vert \V(x)-\V(0)\vert}{\vert x \vert^{1-\frac{3}{p}}}\leq \Vert \V \Vert_{\dot{B}^{1-\frac{3}{p}}_{\infty,\infty}}},$$ and by recalling that the vector field $\V_2$ is defined in the second identity in (\ref{def-V1-V2}), then we obtain:
\[\sup_{\vert x \vert <R} \frac{\vert \V_2(x)\vert}{\vert x \vert^{1-\frac{3}{p}}}\leq \Vert \V \Vert_{\dot{B}^{1-\frac{3}{p}}_{\infty,\infty}}. \]
Thereafter, for all $x \in B_R$ we have:  
\[ \vert \V_2(x)\vert \leq \Vert \V \Vert_{\dot{B}^{1-\frac{3}{p}}_{\infty,\infty}} \, \vert x \vert^{1-\frac{3}{p}} \leq  \Vert \V \Vert_{\dot{B}^{1-\frac{3}{p}}_{\infty,\infty}}  \,R^{1-\frac{3}{p}}. \]
We thus have:
\begin{equation*}
\begin{split}
	I(R) \leq & \, C\,  \frac{1}{R^4} \left(\frac{1}{R^3} \int_{\vert x \vert <R} \vert \V_2(x) \vert^{p}dx \right)^{\frac{4}{p-3}} \leq \Vert \V \Vert_{\dot{B}^{1-\frac{3}{p}}_{\infty,\infty}} \,  \frac{1}{R^4} \left( \left( \frac{1}{R^3} \int_{\vert x \vert <R} dx \right) R^{p(1-\frac{3}{p})}  \right)^{\frac{4}{p-3}}  \\ 
\leq & \, 	C÷,  \Vert \V \Vert_{\dot{B}^{1-\frac{3}{p}}_{\infty,\infty}} \,  \frac{1}{R^4}\left( R^{p-3} \right)^{\frac{4}{p-3}} \leq c \Vert \V \Vert_{\dot{B}^{1-\frac{3}{p}}_{\infty,\infty}} \leq c \Vert \U \Vert_{\dot{B}^{-3/p}_{\infty, \infty}} \leq c \Vert \U \Vert_{L^p}. 
\end{split}
\end{equation*}

 We get back to (\ref{Caccioppoli2}) to write 
	$$ \int_{B_{R/2}} \vert \vec{\nabla} \otimes \U(x) \vert^2 dx \leq C_p  \left( \int_{\mathcal{C}(R/2,R)} \vert \V_1(x) \vert^2 dx \right) \Vert \U \Vert_{L^p}.$$ 
We recall  that by the first identity in  formula (\ref{def-V1-V2}) we have $\V_1=\V$ where $\ds{\V= \frac{1}{-\Delta} (\vec{\nabla} \wedge \U)}$. Moreover,   as we have  $\U \in \dot{H}^{-1}(\Rt)$ then we obtain  $\V \in L^{2}(\Rt)$. Consequently, by  letting	$R \to +\infty$ we have  $\Vert \vec{\nabla} \otimes \U \Vert^{2}_{L^2}=0$. By proceeding as in the end of the proof of Theorem \ref{Th1} we have $\U=0$ and $P=0$.  Theorem \ref{Th2} is  proven. \finpv
	
\subsection{Proof of Corollary \ref{Col}}
Let $\U \in \dot{H}^{1}(\Rt)$ be a weak solution of the (NS) equations. First we remark that as $\U \in \dot{H}^{1}(\Rt)$, by the Sobolev embeddings we have $\U \in L^6(\Rt)$, and consequently, $\U \in L^{3}_{loc}(\Rt)$. Then, by Theorem $X.1,1$ at the  page 658 in \cite{Galdi} we have   $\U \in \mathcal{C}^{\infty}(\Rt)$ and $P \in \mathcal{C}^{\infty}(\Rt)$. 

\medskip 

Now, we  shall prove that $\U \in \dot{H}^{-1}(\Rt)$.  For  $\rho>0$ fixed, we write 
\[\Vert \U \Vert^{2}_{\dot{H}^{-1}}= \int_{\Rt} \frac{1}{\vert \xi \vert^2} \left\vert \widehat{\U} (\xi)\right\vert^2 \, d\xi =  \int_{\vert \xi \vert < \rho} \frac{1}{\vert \xi \vert^2} \left\vert \widehat{\U} (\xi)\right\vert^2 \, d\xi + \int_{\vert \xi \vert \geq \rho} \frac{1}{\vert \xi \vert^2} \left\vert \widehat{\U} (\xi)\right\vert^2 \, d\xi.  
\]
In order to control the first term  on the right-hand side, by the H\"older inequalities (with $1=2/p+2/r$) we write 
\[ \int_{\vert \xi \vert < \rho} \frac{1}{\vert \xi \vert^2} \left\vert \widehat{\U} (\xi)\right\vert^2 \, d\xi \leq \left( \int_{\vert \xi \vert < \rho} \frac{1}{\vert \xi \vert^{p}}\, d\xi  \right)^{2/p} \left( \int_{\vert \xi \vert < \rho}  \left\vert \widehat{\U} (\xi)\right\vert^{r}\, d \xi\right)^{2/r}.\]
We set $p<3$ and  the first term in the right ride converges. Moreover, the hypothesis $\widehat{\U} \in L^{r}_{loc}(\Rt)$, with $r>6$, allows us to conclude that the second term on the right-hand side also converges. 

\medskip 

To control the second term on the right-hand side of the last identity, we just write 
\[  \int_{\vert \xi \vert \geq \rho} \frac{1}{\vert \xi \vert^2} \left\vert \widehat{\U} (\xi)\right\vert^2 \, d\xi = \int_{\vert \xi \vert \geq \rho} \frac{1}{\vert \xi \vert^4} \vert \xi \vert^2 \left\vert \widehat{\U} (\xi)\right\vert^2 \, d\xi \leq \frac{1}{\rho^4} \Vert \U \Vert^{2}_{\dot{H}^1}<+\infty.\]

We this get $\U \in L^6(\Rt) \cap \dot{H}^{-1}(\Rt)$ and  by Theorem \ref{Th2} we have the identities  $\U=0$ and $P=0$.  \finpv 
\begin{appendices}
\section{Appendix}\label{AppendixA}

Let $f \in \mathcal{C}^3(\Rt)$ be a harmonic function. We define $\U= \vec{\nabla} f $ and $\ds{P = - \frac{1}{2} \vert \vec{\nabla} f \vert^2}$, and we will prove that $(\U,P)$ is a solution of the (NS) equations. Indeed, as $\Delta f =0$ we directly have $-\Delta \U=0$. On the other hand, by well-known rules of the vector calculus we have the identity 
\[(\U \cdot \vec{\nabla}) \U= \vec{\nabla} \left( \frac{1}{2}\vert \U \vert^2 \right)+ \left(\vec{\nabla} \wedge \U\right)\wedge \U, \]
 but, as $\U= \vec{\nabla} f $ then we have $\ds{\vec{\nabla} \wedge \U=0}$, and we can write 
 \[ (\U \cdot \vec{\nabla}) \U=  \frac{1}{2} \vec{\nabla}\left(\vert \U \vert^2\right) =-\vec{\nabla} P.\]
 Hence, $(\U,P)$ solves the equation $\ds{-\Delta \U + (\U \cdot \vec{\nabla}) \U +\vec{\nabla} P = 0}$.  Moreover, always as by the fact that $f$ is an harmonic function we have $\ds{div(\U)= div (\vec{\nabla} f)=\Delta f =0}$. 
\end{appendices}

\section*{Availability of data and materials}
Data sharing not applicable to this article as no datasets were generated
or analyzed during the current study.

\vspace{5cm} 

\end{document}